\documentclass[twoside,reqno,a4paper,12pt]{amsart}
\usepackage{amsmath,amsthm}
\usepackage{amsmath,amstext,amsfonts}
\usepackage{hyperref}
\usepackage{graphicx}
\usepackage{graphics}
\usepackage{tikz}
\usepackage{float}
\usepackage{etex}

\theoremstyle{plain}

\newtheorem{theo}{Theorem}

\theoremstyle{remark}

\def \to {\rightarrow}

\newcommand{\trip}[1]{{|\kern -1pt\|#1\|\kern -1pt|}}

\begin{document}

\title[An example on Lyapunov stability and linearization]
{An example on Lyapunov stability and linearization}

\author[Hildebrando M. Rodrigues]
{Hildebrando M. Rodrigues$^{\dag}$}

\thanks{\dag Instituto de Ci\^{e}ncias Matem\'{a}ticas e de
Computa\c{c}ao, Universidade de S\~{a}o Paulo-Campus de S\~{a}o
Carlos, Caixa Postal 668, 13560-970 S\~{a}o Carlos SP, Brazil, e-mail: hmr@icmc.usp.br}

\author[J. Sol\`a-Morales]
{J. Sol\`a-Morales$^{\ddag}$}
\thanks{ \ddag Departament de
Matem\`atiques, Universitat Polit\`ecnica de Catalunya,
Av. Diagonal 647, 08028 Barcelona, Spain, e-mail:
jc.sola-morales@upc.edu}
 \markboth{Joan de Sol\`a-Morales}{
Hildebrando M. Rodrigues}
\thanks{\dag Partially supported by FAPESP Processo 2018/05218-8 }
\thanks{\ddag Partially supported by MINECO grant MTM2017-84214-C2-1-P. Faculty member of the Barcelona Graduate School of Mathematics (BGSMath) and part of the Catalan research group 2017 SGR 01392.}

\maketitle

\begin{abstract} The purpose of this paper is to present an example of a $\mathcal{C}^1$ (in the Fr\'{e}chet sense) discrete dynamical system in a infinite-dimensional separable Hilbert space for which the origin is an exponentially asymptotically stable fixed point, but such that its derivative at the origin has spectral radius larger than unity, and this means that the origin is unstable in the sense of Lyapunov for the linearized system. The possible existence or not of an example of this kind has been an open question until now, to our knowledge. The construction is based on a classical example in Operator Theory due to Kakutani.

\end{abstract}

\noindent 

\noindent 

\vspace{0.5in}

\section{Introduction and description of Kakutani's example}

Let $X$ be a real Banach space and $T:X\to X$ a map such that $T(0)=0$ and $T$ is differentiable in the Fr\'{e}chet sense at $x=0$. If we call $\mathcal{M}=DT(0)$, then we can write $T(x)=\mathcal{M}x+N(x)$ and the nonlinear part $N(x)$ satisfies $\|N(x)\|/\|x\|\to 0$ as $x\to 0$. It is well-known and very easy to prove that if the spectral radius of $\mathcal{M}$ is less than one then the origin is exponentially assimptotically stable, in the sense of Lyapunov, as a fixed point of $\mathcal{M}$ and as a fixed point of $T$. If the spectral radius of $\mathcal{M}$ is larger than one, then the origin is exponentially unstable for the linear system defined by $\mathcal{M}$, but the situation is not so clear concerning the instability for the nonlinear system defined by $T$.

In a recent very interesting paper, entitled {\sl On Nonlinear Stabilization of Linearly Unstable Maps} \cite{GTZ}, the authors, T. Gallay, B. Texier and K. Zumbrun, underline, even in the Abstract, that their results {\sl ``highlight the fundamental open question whether Fr\'{e}chet differentiability is sufficient for linear exponential instability to imply nonlinear exponential instability, at possibly slower rate"}. Later on, in the Introduction section, they say that they {\sl ``are not aware of any example of nonlinear stabilization for a linearly unstable Fr\'{e}chet differentiable dynamical system, nor of any result that would prevent such a phenomenon to occur under minimal and natural assumptions"}. In \cite{GTZ} the above authors give interesting examples of nonlinear stabilization, but only with Gateaux differentiability.

The goal of the present article is to present such an example, with Fr\'{e}chet differentiability.

The main instability result for general linear parts still remains to be, to our knowledge, the instability theorem of Dan Henry (\cite{Henry}, Theorem 5.1.5, p. 105), namely

(D. Henry, 1981): {\sl If $X$ is a Banach space and $T_n$ is a continuous map from a neighbourhood of the origin of $X$ into $X$ with $T_n(0)=0$ $(n=1,2,3,\dots)$ and $\mathcal{M}$ is a continuous linear operator on $X$ with spectral radius greater than one, and
\begin{equation}\label{Henry}
T_n(x)=\mathcal{M}x+O(\|x\|^p) \ \ \ \hbox{ as } x\to 0
\end{equation}
uniformly in $n\ge 1$, for some constant $p>1$, then the origin is unstable. Specifically, there exists a constant $C>0$ and there exists $x_0$ arbitrarily close to $0$ such that, if $x_n=T_n(x_{n-1})$ for $n\ge 1$, then for some $N$ (depending on $x_0$), the sequence $x_0,x_1,\dots,x_N$ is well defined and $\|x_N\|\ge C$.}

Observe that in this statement the nonlinear part $T_n(x)-\mathcal{M}x$ may depend on $n$ as long as the bound $O(\|x\|^p)$ remains uniform on $n$. Henry's result, on discrete non-autonomous dynamical systems, is the main ingredient of his Corollary 5.1.6, were it is applied to semilinear evolution equations, mostly semilinear parabolic partial differential equations, that is the principal goal of his book. As it is said in (\cite{Henry} page 331), his Theorem is inspired in Theorem 2.3 (Ch. 7) of \cite{DK} (p. 291), where the situation is considered in the context of continuous-time dynamical systems defined by ordinary differential equations in Banach spaces, with a bounded linear part. See also the notes on p. 332 of \cite{DK} that refer to a previous proof by M.A. Rutman. This $O(\|x\|^p)$ with $p>1$ bound for the nonlinearity is also a key hypothesis in the instability result of J. Shatah and W. Strauss \cite{SS}.

In \cite{GTZ} an improvement of the previous theorem is also presented, under a new sufficient condition, but only for the case of Hilbert spaces and when the linear part $\mathcal{M}$ is self-adjoint or normal. In that case, condition \eqref{Henry} is slightly weakened to
\begin{equation}\label{weakened}
  \|T(x)-\mathcal{M}x\|\le\alpha(\|x\|)\|x\|\ \ \hbox{ as } x\to 0, \text{ with } \int_0^a\dfrac{\alpha(s)}{s}\ ds<\infty, \text{ for some } a>0.
\end{equation}

We want to say that our work started after attending the stimulating talk {\sl Remarks on Nonlinear Stabilization} of Prof. T. Gallay, where he explained the results of \cite{GTZ}. This happened in the conference {\sl Dynamics of Evolution Equations}, a conference dedicated to the memory of Prof. Jack K. Hale, which took place at the CIRM, in Luminy, in March 2016. We want also to say that we have been interested in linearization in infinite dimensions for a number of years, starting with our work \cite{Ro-Sola}. See a summary of our work in \cite{6Ro-Sola}. And the present paper is not the first time we find an example of a phenomenon in this field that cannot happen in finite dimensions, as it was the example we studied in \cite{3Ro-Sola}.

There is a completely different approach to instability results by linearization, perhaps more close to Lyapunov original ideas, that uses the additional hypothesis of the {\sl spectral gap condition} or the separation of the spectrum of the linear part by a circle around zero of radius one or larger than one (see, \cite{Iooss}, Ch. I, Th. 2, p. 3). This happens automatically, for example, in finite dimensions. Something like our example cannot be expected in a situation like this.

Let us describe the principal features of our example in the following statement, that summarizes our result.

\begin{theo}
In an infinite-dimensional separable real Hilbert space $\mathcal{H}$ there exist $\mathcal{C}^1$ maps (in the sense of Fr\'{e}chet) $T:\mathcal{H}\to\mathcal{H}$ of the form $T(x)=W_\varepsilon x+N(x)$, where $W_\varepsilon$ is a bounded linear operator of the type of a weighted shift, on a given Hilbert basis, and $N(x)$ satisfies that $\|N(x)\|/\|x\|\to 0$ as $x\to 0$ with the property that the spectral radius of $W_\varepsilon$ is larger than one but the fixed point $x=0$ is exponentially asymptotically stable. In fact, these maps will depend on a real sequence $(\varepsilon_m)$ with $\varepsilon_m\searrow 0$, to be described later, and one will have
\begin{equation}\label{finebound}
  \left(\dfrac{1}{-\log\|x\|}\right)^{c_2}<\dfrac{\|N(x)\|}{\|x\|}<4\,\left(\dfrac{1}{-\log\|x\|}\right)^{c_1}
\end{equation}
as $x\to 0$ for some $0<c_1<c_2$ that can be chosen arbitrarily if one chooses accordingly the sequence $(\varepsilon_m)$ .
\end{theo}

Observe that if one chooses $c_1>1$ in \eqref{finebound}, then condition \eqref{weakened} holds. This means that \eqref{weakened} is a sufficient condition for the case of self-adjoint or normal linear parts, but not for more general linear operators like our weighted shifts $W_\varepsilon$.

Since there are slight differences in the literature, it is worth to say that we are using the following definition of exponential asymptotic stability of $x=0$ as a fixed point of $T$ (see T. Yoshizawa, \cite{Yoshizawa}, p. 48, Definition 7.6): There exists a $\lambda<1$ such that for any $\varepsilon>0$ there corresponds a $\delta(\varepsilon)>0$ such that $\|x_0\|<\delta(\varepsilon)$ implies $\|T^n(x_0)\|\le\varepsilon\lambda^n$, for all $n\ge 0$.

The rest of this paper is organized as follows. In the next paragraphs we describe an important construction due to S. Kakutani in the context of Operator Theory, which is the basis of our example. In fact, our example can be seen as a nonlinear version of Kakutani's. In Section 2 we construct our nonlinear map $T$, we study its regularity and we obtain the bounds and properties that will be needed in the sequel. In Section 3 we prove the stability and the exponential asymptotic stability of its fixed point $x=0$.

Let's start with the description of Kakutani's construction. In an infinite-dimensional separable Hilbert space $\mathcal{H}$ with a Hilbert basis $(e_n)_{n\ge 1}$ a {\sl weighted shift} is a bounded linear operator $W\in\mathcal{L}(\mathcal{H})$ that is defined by $We_n=\alpha_ne_{n+1}$, for a bounded sequence $(\alpha_n)_{n\ge 1}$. It is easy to prove that $\|W\|=\sup_n\{|\alpha_n|\}$ and that $\|W^k\|=\sup_{n}\{|\alpha_n\alpha_{n+1}\cdots\alpha_{n+k-1}|\}$. In matrix form it would be
$$W=\begin{pmatrix}
0&0&0&0&\cdots\\
\alpha_1&0&0&0&\cdots\\
0&\alpha_2&0&0&\cdots\\
0&0&\alpha_3&0&\cdots\\
\vdots&\vdots&\vdots&\vdots&\vdots\\
\end{pmatrix}
,\ \
W^2=\begin{pmatrix}
0&0&0&0&\cdots\\
0&0&0&0&\cdots\\
\alpha_1\alpha_2&0&0&0&\cdots\\
0&\alpha_2\alpha_3&0&0&\cdots\\
\vdots&\vdots&\vdots&\vdots&\vdots\\
\end{pmatrix}
,\ \ \cdots
$$

The example of Kakutani is a weighted shift operator with positive spectral radius that can be arbitrarily approximated in the operator norm by nilpotent operators, that are operators with spectral radius equal to zero. Following \cite{Rickart} (p. 282) and \cite{Halmos} (p. 248), Kakutani's example starts with a strictly decreasing sequence of positive real numbers $(\varepsilon_m)$ ($m\ge 1)$ converging to $0$, and with another property to be stated later.

In general, every positive integer $n\ge1$ can be written, and in a unique way, as the product of a power of $2$ with an odd integer, that is $n=2^k(2\ell+1)$, with $k,\ell\ge 0$ nonnegative integers. Therefore, we can consider the integer-valued function $k=k(n)$.  A second sequence $(\alpha_n)$ ($n\ge 1)$ is defined now by $\alpha_n=\varepsilon_{1+k(n)}$. It is clear that a member $\varepsilon_\mu$ of the $(\varepsilon_m)$ sequence will appear infinitely many times in the sequence $(\alpha_n)$, and precisely in the places $n=2^{(\mu-1)}(2\ell+1)$, for $\ell=0,1,\dots$. Therefore, $\varepsilon_\mu$ will appear periodically in the sequence $(\alpha_n)$, with a period $2^\mu$, and it will appear for the first time at the position $n=2^{(\mu-1)}$.

The sequence $(\alpha_n)$ for $n=1,\dots, 31\dots$, for example, looks like this:
\begin{equation*}
  \varepsilon_1, \varepsilon_2, \varepsilon_1, \varepsilon_3, \varepsilon_1, \varepsilon_2, \varepsilon_1, \varepsilon_4, \varepsilon_1, \varepsilon_2, \varepsilon_1, \varepsilon_3, \varepsilon_1, \varepsilon_2, \varepsilon_1, \varepsilon_5,
\varepsilon_1, \varepsilon_2, \varepsilon_1, \varepsilon_3, \varepsilon_1, \varepsilon_2, \varepsilon_1,
\end{equation*}
\begin{equation*}
  \varepsilon_4, \varepsilon_1, \varepsilon_2, \varepsilon_1, \varepsilon_3, \varepsilon_1, \varepsilon_2, \varepsilon_1, \dots
\end{equation*}

Let now $\mathcal{H}$ and $(e_n)$ be as above and define the weighted shift $W_\varepsilon$ by $W_\varepsilon e_n=\alpha_ne_{n+1}$ with the sequence $(\alpha_n)$ just defined. It is clear that $\|W_\varepsilon\|=\varepsilon_1$. Define also the weighted shifts $L_m$ by the same formula as $W_\varepsilon$ but with all the $\varepsilon_k$ replaced by $0$ except $k=m$. The sequence of weights for $L_3$, for example, would be
$$
0, 0, 0, \varepsilon_3, 0, 0, 0, 0, 0, 0, 0, \varepsilon_3, 0, 0, 0, 0, 0, 0, 0, \varepsilon_3, 0, 0, 0, 0, 0, 0, 0, \varepsilon_3, \cdots.$$
It is clear that $\|L_m\|\to 0$ as $m\to\infty$ and that $W_\varepsilon-L_m$ is a weighted shift with weights equal to zero each $2^m$ places, starting at the place $2^{(m-1)}$. Therefore, $W_\varepsilon-L_m$ is nilpotent of index $2^m$. In conclusion, $W_\varepsilon$ can be approximated in the operator norm by a sequence $W_\varepsilon-L_m$ of nilpotent operators. Observe, for later use, that
\begin{equation}\label{later}
  \|W_\varepsilon-L_m\|\le\|W_\varepsilon\|, \ \hbox{ for all }\  m,
\end{equation}
since $W_\varepsilon-L_m$ is also a weighted shift.

Let us see now that for a good choice of the sequence $(\varepsilon_m)$ the spectral radius $\rho(W_\varepsilon)$ of $W_\varepsilon$ turns out to be positive. Observe that $\|W_\varepsilon^n\|=\alpha_1\alpha_2\cdots\alpha_n$ and that for $p\ge 1$ and $n=2^p-1$
$$\|W_\varepsilon^n\|=\|W_\varepsilon^{2^p-1}\|=\varepsilon_{1}^{2^{(p-1)}}\varepsilon_2^{2^{(p-2)}}\cdots\varepsilon_{p-2}^4\varepsilon_{p-1}^2\varepsilon_p^1=\prod_{q=1}^{p}\varepsilon_q^{2^{(p-q)}}.$$
Hence,
$$\log \left(\|W_\varepsilon^n\|^{\frac{1}{n}}\right)=\dfrac{1}{2^p-1}\sum_{q=1}^p2^{(p-q)}\log\varepsilon_q=\dfrac{2^p}{2^p-1}\sum_{q=1}^p \dfrac{\log\varepsilon_q}{2^q}.$$

Several choices of the $(\varepsilon_m)$ sequence will make this series to be convergent as $p\to\infty$, and $\rho(W_\varepsilon)>0$. This is Kakutani's example: an operator $W_\varepsilon$ with positive spectral radius that can be arbitrarily approximated in the operator norm by nilpotent operators, $W_\varepsilon-L_m$ which in particular have spectral radius equal to zero. We are going to use these ideas and this notation to build our nonlinear map.

\section{Construction of the map}

Following the previous notation, let us take now $(\varepsilon_m)$ as $\varepsilon_m=M/K^{m-1}$ with $M>K>1$. Then, $\|W_\varepsilon\|=M>1$ and it is easy to calculate that $\rho(W_\varepsilon)=M/K>1$, since $\sum_{q=1}^\infty (1-q)2^{-q}=-1$. Therefore, the linear dynamical system $x_{n+1}=W_\varepsilon x_n$ in $\mathcal{H}$ has $x=0$ as an unstable fixed point.

To fix ideas, a choice like $M=5$ and $K=3$, for example, would be appropriate for all what is said in the rest of the paper. But we continue with the use of the generic letters $M$ and $K$ with the assumption that $M>K>1$.

We are going to construct now a nonlinear map $N:\mathcal{H}\to\mathcal{H}$ of class $\mathcal{C}^1$ such that $N(x)=o(\|x\|)$ as $x\to 0$ and such that if we define $T:=W_\varepsilon+N$ then $x=0$ is a fixed point of $T$ that is exponentially asymptotically stable. This map will have the form $N(x)=\widetilde{N}(\|x\|)x$ for a $\mathcal{C}^1$ map $(0,\infty)\ni t\mapsto \widetilde{N}(t)\in\mathcal{L}(\mathcal{H})$, and $N(0)=0$.

Given two real numbers $a<b$ let $\varphi_1(a,b;t)$ be a real function of $t\in[a,b]$ of class $\mathcal{C}^1$ in $t$ such that $\varphi_1(a,b;a)=\partial_t\varphi_1(a,b;a)=0$, $\varphi_1(a,b;b)=1$, $\partial_t\varphi_1(a,b;b)=0$ and $0\le\varphi_1(a,b;t)\le 1$ for all $t\in[a,b]$. We can suppose also that $0\le\partial_t\varphi_1(a,b;t)\le 2/(b-a)$. Let us define also $\varphi_2=1-\varphi_1$, that reverses the values at the limits $a$ and $b$.

\begin{figure}[h]
\begin{tikzpicture}[xscale=10,yscale=1/3]
\draw[->,black, thin, domain=0:1.3] plot (\x,0);
\node at (1.35,-.05) {$t$};
\draw[black,thin] ( .23,0.4) -- ( .23,-0.4) node[below,black]{$M^{-2^{(k+3)}}$};
\draw[black,thin] ( .38,0.4) -- ( .38,-0.4) node[below,black]{$M^{-2^{(k+2)}}$};
\draw[black,thin] ( .59,0.4) -- ( .59,-0.4) node[below,black]{$M^{-2^{(k+1)}}$};
\draw[black,thin] ( .82,0.4) -- ( .82,-0.4) node[below,black]{$M^{-2^{k}}$};
\draw[black,thin] (1.08,0.4) -- (1.08,-0.4) node[below,black]{$M^{-2^{(k-1)}}$};
\draw[black, thick, domain=0:0.23] plot (\x,.08+0);
\draw[black, thick, domain=0.23:0.38] plot (\x,{.08  +(5*((\x-0.23)/(.15))^2-10*(\x-0.23)^2*(\x-0.38)/(.15^3)});
\draw[black, thick, domain=0.38:0.82] plot (\x,.08+5);
\draw[black, thick, domain=0.82:1.08] plot (\x,{.08+5-(5*((\x-0.82)/(.26))^2-10*(\x-0.82)^2*(\x-1.08)/(.26^3)});
\draw[black, thick, domain=1.08:1.22] plot (\x,.08+0);
\node at (.595,8.4) {{\sl $\widetilde{N}_k=-L_k$}};
\node at (0.12,4.2) {{$\widetilde{N}_k=0$}};
\node at (1.14,4.2) {{$\widetilde{N}_k=0$}};
\draw[->,black,thin] (0.69,7.5) -- (0.76,6.5);
\draw[->,black,thin] (0.52,7.5) -- (0.45,6.5);
\draw[->,black,thin] (.12,3.2) -- (.12,1.2);
\draw[->,black,thin] (1.14,3.2) -- (1.14,1.2);
\end{tikzpicture}
\caption{Sketch of the graph of $\widetilde{N}_k(t)$}
\label{graph}
\end{figure}
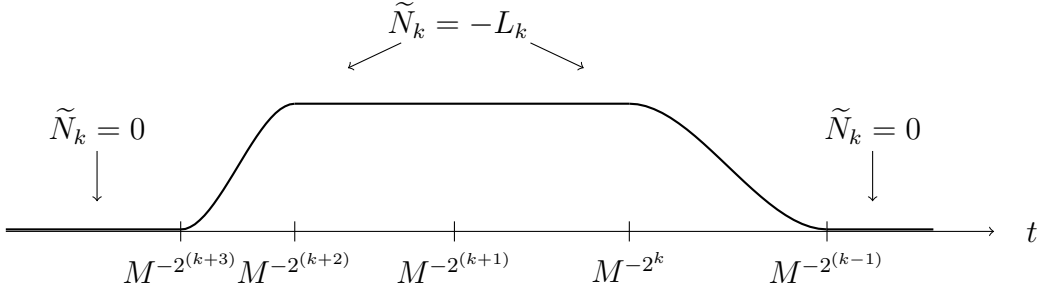

Let us define now, for $t> 0$ and $k\ge 1$
\begin{equation}\label{Nk}
  \widetilde{N}_k(t)=
  \begin{cases}
  0,\text{ if } t<M^{-2^{(k+3)}}\\
  -\varphi_1(M^{-2^{(k+3)}},M^{-2^{(k+2)}};t)L_k, \text{ if } t\in[M^{-2^{(k+3)}},M^{-2^{(k+2)}})\\
  -                                              L_k, \text{ if } t\in[M^{-2^{(k+2)}},M^{-2^{ k   }})\\
  -\varphi_2(M^{-2^{ k   }},M^{-2^{(k-1)}};t)L_k, \text{ if } t\in[M^{-2^{ k   }},M^{-2^{(k-1)}})\\
   0,\text{ if } t\ge M^{-2^{(k-1)}}.

  \end{cases}
\end{equation}
 See Fig. \ref{graph} for a sketch of the graph of $\widetilde{N}_k(t)$. It is clear that $\widetilde{N}_k:(0,\infty)\to\mathcal{L}(\mathcal{H})$ is a $\mathcal{C}^1$ function such that for $M^{-2^{(k+2)}}\le t<M^{-2^{ k   }}$ one has that $(W_\varepsilon+\widetilde{N}_k(t))$ is a nilpotent linear operator of index $2^k$. To describe better the situation, we need to introduce a new property, that could be called the {\em nilpotency of a set of operators} $\Omega\subset\mathcal{L}(\mathcal{H})$.

For $k\ge 1$ we define $\Omega_k\subset\mathcal{L}(\mathcal{H})$ as the class of weighted shifts defined by a sequence of weights $(\alpha_n)$ such that $\alpha_n=0$ for $n=2^{(k-1)}(2\ell+1)$ and $\ell=0,1,\dots$ In other words, the class $\Omega_k$ consists of the weighted shifts that have weights $\alpha_n$ equal to zero at least in the place $n=2^{(k-1)}$ and from this place onwards each $2^k$ positions. Observe that these classes are not disjoint and that, according to \eqref{Nk}, for example, if
$\|x\|\in[M^{-2^{(k+1)}},M^{-2^{ k   }})$ then
\begin{equation}\label{intersection}
  \left(W_\varepsilon+\widetilde{N}_{k+1}(\|x\|)+\widetilde{N}_k(\|x\|)+\widetilde{N}_{k-1}(\|x\|)+\widetilde{N}_{k-2}(\|x\|)\right)
  \in\Omega_k\cap\Omega_{k-1},
\end{equation}
for $k\ge 3$.

We remark now the nilpotency of the set of operators $\Omega_k$, in the sense that
\begin{equation}\label{nil}
  \hbox{ if  }\ \ W_1,W_2,\dots W_{2^k}\in \Omega_k, \ \ \hbox{ then } \ \ W_1\cdot W_2\cdots W_{2^k}=0,
\end{equation}
as it is easily proved: (we call $(\alpha_i)_j$ the $i$-th weight of $W_j$) for all $i\ge 1$, $\ W_1\cdot W_2\cdots W_{2^k}(e_i)=(\alpha_i)_{2^k}(\alpha_{i+1})_{2^k-1}\cdots(\alpha_{i+2^k-1})_1e_{i+2^k}$, and this last product must contain at least one factor that is zero. This is a generalization of the observation already made above that since $W_\varepsilon-L_m\in \Omega_m$ then $W_\varepsilon-L_m$ is nilpotent of index $2^m$.

 Let us now define the maps $\widetilde{N}(t):=\sum_{k=1}^\infty \widetilde{N}_k(t)$, $N_k(x):=\widetilde{N}_k(\|x\|)x$ and $N(x):=\widetilde{N}(\|x\|)x$, or, in other words
 \begin{equation}\label{N}
   N(x)=\sum_{k=1}^\infty \widetilde{N}_k(\|x\|)x.
 \end{equation}
Because of \eqref{later} and the fact that $0\le\varphi_1,\varphi_2\le 1$ we see that
\begin{equation}\label{bound}
  \|T(x)\|=\|W_\varepsilon x+N(x)\|\le M\|x\|
\end{equation}
for all $x\in\mathcal{H}$, because $W_\varepsilon+\sum_{k=1}^\infty \widetilde{N}_k(\|x\|)$, for fixed $\|x\|$, is a weighted shift operator whose weights never exceed $M$. We also see that for every $x\in\mathcal{H}$, $x\ne 0$, the expression \eqref{N} is the sum of a maximum of four nontrivial or active terms $\widetilde{N}_k(\|x\|)x$, that are also consecutive, namely, if
$\|x\|\in[M^{-2^{(k+1)}},M^{-2^{ k   }})$ then
$$
  N(x)=\widetilde{N}_{k+1}(\|x\|)x+\widetilde{N}_k(\|x\|)x+\widetilde{N}_{k-1}(\|x\|)x+\widetilde{N}_{k-2}(\|x\|)x,
$$
for $k\ge 3$, and there exists a small neighbourhood of $x$ where
the expression \eqref{N} is the sum of a maximum of five active and consecutive terms $\widetilde{N}_k$. This last observation implies that $N$ is of class $\mathcal{C}^1$ at least in $\mathcal{H}\setminus\{0\}$.

To prove that it is also $\mathcal{C}^1$ up to $x=0$ we start by showing that $DN(0)=0$. Let $x\ne 0$ and suppose that $\|x\|\in[M^{-2^{(k+1)}},M^{-2^{ k   }})$, with $k\ge 3$. Then,
$$\dfrac{\|N(x)\|}{\|x\|}\le \dfrac{\|L_{k+1}x\|+\|L_{k}x\|+\|L_{k-1}x\|+\|L_{k-2}x\|}{\|x\|}$$ $$\le \varepsilon_{k+1}+\varepsilon_{k}+\varepsilon_{k-1}+\varepsilon_{k-2}
<4\,\varepsilon_{k-2}$$
and this tends to $0$ as $x\to 0$, because $k\to\infty$ as $x\to 0$. This proves that $N(x)=o(\|x\|)$ and therefore it is Fr\'{e}chet differentiable at $x=0$ with $DN(0)=0$.

Let us make more explicit the dependence of $\varepsilon_k$ on $x$ to prove \eqref{finebound}.
First of all, let us see that if $R\in[M^{-2^{(k+1)}},M^{-2^{ k   }})$, and, as before, $(e_n)$ is the Hilbert basis, then $N(R\,e_{2^{(k-1)}})=-\varepsilon_k R\,e_{2^{(k-1)}+1}$. This is so because for $x=Re_{2^{(k-1)}}$ only the terms $\widetilde{N}_{k+1}, \widetilde{N}_{k}, \widetilde{N}_{k-1}$, and $\widetilde{N}_{k-2}$ are active in \eqref{N}, and when applied to $R\,e_{2^{(k-1)}}$ only the one of index $k$ is nonzero. Exactly, it gives $\widetilde{N}_k(R)R\,e_{2^{(k-1)}}=-L_kR\,e_{2^{(k-1)}}=-\alpha_{2^{(k-1)}}R\,e_{2^{(k-1)}+1}=-\varepsilon_k R\,e_{2^{(k-1)}+1}$.

Therefore, the previous inequality can be improved to the following: If $R\in[M^{-2^{(k+1)}},M^{-2^{ k   }})$, then
\begin{equation}\label{finebounds}
  \varepsilon_{k}\le\sup_{\|x\|=R}\dfrac{\|N(x)\|}{\|x\|}
  \le 4\,\varepsilon_{k-2}.
\end{equation}

We have $\varepsilon_k=M/K^{k-1}$, and therefore $(-\log(\varepsilon_k))/k\to\log K$ as $k\to\infty$. The same calculation gives, $(-\log(\varepsilon_{k-2}))/k\to\log K$ as $k\to\infty$. In addition, $M^{-2^{(k+1)}}\le\|x\|<M^{-2^{ k   }}$, and we see that $(\log(-\log\|x\|))/k\to\log 2$. Then, there exist positive numbers $0<c_1<c_2$ such that $c_1(\log(-\log\|x\|))\le -\log(\varepsilon_{k-2})<
-\log(\varepsilon_k)\le c_2(\log(-\log\|x\|))$ and
$$\left(\dfrac{1}{-\log\|x\|}\right)^{c_2}<\varepsilon_k<\varepsilon_{k-2}<\left(\dfrac{1}{-\log\|x\|}\right)^{c_1},$$
for $k$ sufficiently large, or $\|x\|$ sufficiently small, and the bound \eqref{finebound} is satisfied.

More precisely, $c_1$ and $c_2$ need only to satisfy $0<c_1<\log K/\log 2<c_2$, and therefore they are somehow arbitrary since we are allowed to change the value of $K$ as long as $1<K<M$. For example, in the case $M=5$ and $K=3$ mentioned above, since $\log 3 >\log 2$ then $c_1$ can be chosen grater than $1$, and \eqref{weakened} holds.

Let us see now that $DN(x)\to 0$ as $x\to 0$. For this we calculate a bound for $DN_k(x)$. In general, the Fr\'{e}chet derivative of a function of the form $F(x):=\phi(\|x\|)Ax$, for a smooth real function $\phi(t)$ and a bounded linear operator $A$, at $x\ne 0$ is
$$[DF(x)](y)=\phi'(\|x\|)\langle\dfrac{x}{\|x\|},y\rangle Ax+\phi(\|x\|)Ay,$$
and therefore
$$\|DF(x)\|\le|\phi'(\|x\|)|\,\|A\|\,\|x\|+|\phi(\|x\|)|\,\|A\|.$$

Consequently,

\begin{equation}\label{NDNk}
  \|DN_k(x)\|
  \begin{cases}
  =0,\text{ if } \|x\|<M^{-2^{(k+3)}}\\
  \le\dfrac{2}{M^{-2^{(k+2)}}-M^{-2^{(k+3)}}}\dfrac{M}{K^{k-1}}\|x\|+\dfrac{M}{K^{k-1}}, \text{ if } \|x\|\in[M^{-2^{(k+3)}},M^{-2^{(k+2)}})\\
  \le\dfrac{M}{K^{k-1}}, \text{ if } \|x\|\in[M^{-2^{(k+2)}},M^{-2^{ k   }})\\
   \le\dfrac{2}{M^{-2^{(k-1)}}-M^{-2^{k}}}\dfrac{M}{K^{k-1}}\|x\|+\dfrac{M}{K^{k-1}}, \text{ if } \|x\|\in[M^{-2^{ k   }},M^{-2^{(k-1)}})\\
   =0,\text{ if } \|x\|\ge M^{-2^{(k-1)}},
  \end{cases}
\end{equation}
also
\begin{equation}\label{N2DNk}
  \|DN_k(x)\|
  \begin{cases}
  =0,\text{ if } \|x\|<M^{-2^{(k+3)}}\\
  \le\dfrac{2}{1-M^{-2^{(k+2)}}}\dfrac{M}{K^{k-1}}+\dfrac{M}{K^{k-1}}, \text{ if } \|x\|\in[M^{-2^{(k+3)}},M^{-2^{(k+2)}})\\
  \le\dfrac{M}{K^{k-1}}, \text{ if } \|x\|\in[M^{-2^{(k+2)}},M^{-2^{ k   }})\\
  \le\dfrac{2}{1-M^{-2^{(k-1)}}}\dfrac{M}{K^{k-1}}+\dfrac{M}{K^{k-1}}, \text{ if } \|x\|\in[M^{-2^{ k   }},M^{-2^{(k-1)}})\\
   =0,\text{ if } \|x\|\ge M^{-2^{(k-1)}},
  \end{cases}
\end{equation}
and finally $$\|DN_k(x)\|\le \dfrac{2}{1-M^{-2^{(k-1)}}}\dfrac{M}{K^{k-1}}+\dfrac{M}{K^{k-1}}$$ if $N_k$ is active in $x$. This implies that $\|DN(x)\|\to 0$ as $x\to 0$, because the four indices $k$ that are active in the expression of $N(x)$ tend to $\infty$ as $x\to 0$. This concludes the proof that $N(x)$ is of class $\mathcal{C}^1(\mathcal{H})$.

\section{The stability properties}

Let us proceed now with the proof that $x=0$ as a fixed point of $T=W_\varepsilon+N$ is stable in the sense of Lyapunov. The idea is to take advantage of the fact that a trajectory $x_{n+1}=T(x_n)$ cannot remain in $\|x\|\in[M^{-2^{(k+1)}},M^{-2^k})$ for $2^k$ consecutive indices $n$, because in this range the map $T$ behaves like a nilpotent linear operator of index (at most) $2^{k}$.

Let us come back for one moment to the case $M=5$ and $K=3$ mentioned before. Again to fix ideas, let us suppose that $k=7$, so our initial condition $x_0$ has $\|x_0\|<5^{-2^8}=5^{-256}$. We want to show that for all $n\ge 0$ we will have that $\|x_n\|<5^{-2^7}=5^{-128}$. The idea is that to reach $5^{-128}$, or more, starting below $5^{-256}$, it will take at least $2^7=128$ consecutive steps between $5^{-256}$ and $5^{-128}$, because from $x_n$ to $x_{n+1}$ the norm gets multiplied at most by $5$. And this is impossible, since in between of these bounds, the map $T$ behaves like a nilpotent operator of index at most $128$.

Let us make a precise statement: let us prove that
\begin{equation}\label{stable}
  \hbox{if } \|x_0\|<M^{-2^{(k+1)}} \hbox{ with } k\ge 1 \hbox{ then } \|x_n\|< M^{-2^k} \hbox{ for all } n\ge 0.
\end{equation}
 This implies Lyapunov stability.

Let's argue by contradiction: suppose that for some $n_0$ one has had $\|T^n(x_0)\|< M^{-2^k}$ for all $0\le n<n_0$, but $\|T^{n_0}(x_0)\|\ge M^{-2^k}$. Since  $\|T(x)\|\le M\|x\|$ for all $x$, as we said in \eqref{bound}, we deduce that
$M^{-2^{k}}\le\|T^{n_0}(x_0)\|=\|T^{i}T^{(n_0-i)}(x_0)\|\le M^{i}\,\|T^{(n_0-i)}(x_0)\|$ and
\begin{equation}\label{backn}
\|T^{(n_0-i)}(x_0)\|\ge M^{-2^{k}}M^{-i}, \hbox{ for all } i=0,1,\dots n_0.
\end{equation}
Taking $i=n_0$, \eqref{backn} implies $\|x_0\|\ge M^{-2^k -n_0}$, and since we know that $\|x_0\|< M^{-2^{(k+1)}}$ we deduce that $M^{-2^{k}-n_0}<M^{-2^{(k+1)}}$ and then $n_0>2^{k}$. This means that we can take $i=1,2,\dots 2^{k}$ in \eqref{backn} and obtain that $M^{-2^{k}}>\|T^{(n_0-i)}(x_0)\|\ge M^{-2^{(k+1)}}$ for $i=1,2\dots 2^{k}$.  Therefore, $N_k(x)=-L_kx$ for $x\in\{T^{(n_0-i)}(x_0)|i=1,2\dots 2^{k}\}$, and this would be a contradiction with \eqref{nil}, since that would mean that the $2^k$ operators
$
(W_\varepsilon+\widetilde{N}(\|T^{(n_0-2^k)}(x_0)\|)),
(W_\varepsilon+\widetilde{N}(\|T^{(n_0-2^k+1)}(x_0)\|)),
\dots
(W_\varepsilon+\widetilde{N}(\|T^{(n_0-1)}(x_0)\|))
$ would all of them belong to the same nilpotent set $\Omega_k$, and its product would be $0$. This contradicts that $\|T^{n_0}(x_0)\|\ge M^{-2^k}$.

Let us say again in words the spirit of the last paragraph: the gap between $M^{-2^{k}}$ and $M^{-2^{(k+1)}}$ is so large, relatively to $M$, that if a trajectory starts below $M^{-2^{(k+1)}}$ and at some point arrives above $M^{-2^{k}}$ it must have spent a long time between these two numbers, and this time is large enough to arrive to the index of nilpotency of $W_\varepsilon-L_k$.

Let us write this last result in a different way. Given $x_0$  with $0<\|x_0\|<M^{-4}$ it will exist a $k\ge 1$ such that $M^{-2^{(k+2)}}\le\|x\|<M^{-2^{(k+1)}}$. According to the previous proof, for all $n\ge 0$ one will have that $\|x_n\|<M^{-2^k}$. Therefore,
\begin{equation}\label{stable2}
  \hbox{ if } 0<\|x_0\|<M^{-4},\ \ \hbox{ then } \|T^nx_0\|<\|x_0\|^{1/4}\ \ \hbox{ for all } n\ge 0,
\end{equation}
which is a more quantitative expression of the Lyapunov stability of $x=0$.

To show that $x=0$ is also exponentially asymptotically stable one uses that the functions $\widetilde{N}_k(t)$ defined in \eqref{Nk} are constantly equal to $(-L_k)$ on two consecutive intervals, namely $[M^{-2^{(k+2)}},M^{-2^{(k+1)}})$ and $[M^{-2^{(k+1)}} ,M^{-2^{k}})$.

Given $x_0$ with $\|x_0\|<M^{-1/4}$, a $k(0)\ge 1$ such that $\|x_0\|\in [M^{-2^{(k(0)+2)}},M^{-2^{(k(0)+1)}})$ will necessarily exist. Then, the sequence of $2^{k(0)}$ terms $\|x_0\|$, $\|T(x_0)\|$, $\|T^2(x_0)\|,\dots,$  $\|T^{2^{k(0)}-1}(x_0)\|$ cannot remain in the union of the consecutive intervals $[M^{-2^{(k(0)+2)}},M^{-2^{(k(0)+1)}})$ and $[M^{-2^{(k(0)+1)}} ,M^{-2^{k(0)}})$ without becoming zero. This sequence cannot leave this double interval to the right, because of \eqref{stable}. Therefore it has to leave the interval to the left, and this proves that necessarily an $n(1)$ with $0<n(1)<2^{k(0)}-1$ will exist such that $\|T^{n(1)}x_0\|<M^{-2^{(k(0)+2)}}$.

Let us proceed now inductively with the $k(i)$ and $n(i)$. Observe that if $x_{n(i)}=T^{n(i)}(x_0)\ne 0$ then a $k(i)>k(i-1)$ will necessarily exist such that $M^{-2^{(k(i)+2)}}\le\|x_{n(i)}\|<
M^{-2^{{(k(i)+1)}}}$.  The same argument as above proves that a number $n(i+1)$ will exist such that $n(i)<n(i+1)<n(i)+2^{k(i)}-1$ and $\|T^{n(i+1)}(x_0)\|<M^{-2^{(k(i)+2)}}$.

Along the sequence $(T^n(x_0))_{n\ge 0}$, as long as it does not vanish, the previous inductive step can be repeated, and each term of the sequence $(\|T^n(x_0)\|)_{n\ge 0}$ is bounded from above by the corresponding term of the following sequence:
$$M^{-2^{(k(0)+1)}}, \underbrace{M^{-2^{k(0)}},M^{-2^{k(0)}},\dots M^{-2^{k(0)}}}_{\hbox{ maximum of } 2^{k(0)}-2 \hbox{ terms}}\ \ ,
M^{-2^{(k(1)+1)}}, \underbrace{M^{-2^{k(1)}},M^{-2^{k(1)}},\dots M^{-2^{k(1)}}}_{\hbox{ maximum of } 2^{k(1)}-2 \hbox{ terms}}\ \ ,$$ $$
M^{-2^{(k(2)+1)}}, \underbrace{M^{-2^{k(2)}},M^{-2^{k(2)}},\dots M^{-2^{k(2)}}}_{\hbox{ maximum of } 2^{k(2)}-2 \hbox{ terms}}\ \ ,\ \ \dots
$$
and the terms in this last sequence, are bounded, perhaps not very sharply, but bounded in any case, by the next one
$$\underbrace{M^{-2^{k(0)}},M^{-2^{k(0)}},\dots M^{-2^{k(0)}}}_{\hbox{ exactly } 2^{k(0)}-1 \hbox{ terms}},
\underbrace{M^{-2^{(k(0)+1)}},M^{-2^{(k(0)+1)}},\dots M^{-2^{(k(0)+1)}}}_{\hbox{ exactly } 2^{(k(0)+1)}-1 \hbox{ terms}}\ \ ,$$ $$
\underbrace{M^{-2^{(k(0)+2)}},M^{-2^{(k(0)+2)}},\dots M^{-2^{(k(0)+2)}}}_{\hbox{ exactly } 2^{(k(0)+2)}-1 \hbox{ terms}},\ \ \dots$$
or by the simpler one
$$\underbrace{M^{-2^{k(0)}},M^{-2^{k(0)}},\dots M^{-2^{k(0)}}}_{\hbox{ exactly } 2^{k(0)} \hbox{ terms}},
\underbrace{M^{-2^{(k(0)+1)}},M^{-2^{(k(0)+1)}},\dots M^{-2^{(k(0)+1)}}}_{\hbox{ exactly } 2^{(k(0)+1)} \hbox{ terms}}\ \ ,$$ $$
\underbrace{M^{-2^{(k(0)+2)}},M^{-2^{(k(0)+2)}},\dots M^{-2^{(k(0)+2)}}}_{\hbox{ exactly } 2^{(k(0)+2)} \hbox{ terms}}.\ \ \dots$$
This last sequence has the property that it bounds the sequence $(\|T^n(x_0)\|)_{n\ge 0}$, even if this sequence vanishes after some term. One needs only that
$\|x_0\|\in [M^{-2^{(k(0)+2)}},M^{-2^{(k(0)+1)}})$.

To analyse the last sequence, it is more convenient to analyse only the sequence of the exponents, without the minus sign, and writing $k$ instead of $k(0)$, that is
$$(a_n)_{n\ge 0}:=(\underbrace{{2^{k}},{2^{k}},\dots {2^{k}}}_{ 2^k \hbox{ terms}},
\underbrace{{2^{(k+1)}},{2^{(k+1)}},\dots {2^{(k+1)}}}_{2^{(k+1)} \hbox{ terms}},
\underbrace{{2^{(k+2)}},{2^{(k+2)}},\dots {2^{(k+2)}}}_{2^{(k+2)} \hbox{ terms}},\ \ \dots).$$
We add now the lengths of these blocks and use that $ 2^k+2^{k+1}+2^{k+2}+\dots +2^{k+j}=2^{k+j+1}-2^k,$
and this implies the following values for the final terms of each block: if $n=2^{k+j+1}-2^k$ then $a_n=2^{k+j}$, or in terms of $n$,
$a_n=(n+2^k)/2$. The map $n\mapsto (n+2^k)/2$ is monotone increasing in $n$, coincides with $n\mapsto a_n$ at the end of each block, and $a_n$ is constant on these blocks, therefore we conclude that $a_n\ge (n+2^k)/2$ for all $n\ge 0$.

Therefore, $\|T^n(x_0)\|\le M^{-a_n}\le M^{-(n+2^{k(0)})/2)}$ for all $n\ge 0$, that we can also write as $\|T^n(x_0)\|\le \left(M^{-1/2}\right)^n M^{-2^{(k(0)-1)}}$ or $\|T^n(x_0)\|\le \left(M^{-1/2}\right)^n \left(M^{-2^{(k(0)+2)}}\right)^{1/8}$,
and then
\begin{equation}\label{exp}
  \|T^n(x_0)\|\le \left(M^{-1/2}\right)^n \|x_0\|^{1/8}\ \ \hbox{ for all }\ \ n\ge 0,
\end{equation}
that proves the exponential asymptotic stability, since $M>1$.

\end{document}